\input vanilla.sty
\scaletype{\magstep1}
\scalelinespacing{\magstep1}
\def\bull{\vrule height .9ex width .8ex depth -.1ex}
\def\vol{\text{Vol }}
\def\vr{\text{vr}}
\def\ovr{\text{vr}^*}
\def\dim{\text{dim }}


\pageno=1

\title Factorization theorems for quasi-normed spaces
\endtitle

\author N.J.  Kalton and Sik-Chung Tam\footnote{The research
of both authors was supported by NSF-grant DMS-9201357}\\
Department of Mathematics\\ University of
Missouri-Columbia\\ Columbia, Mo. 65211\\ \endauthor

\vskip1truecm

\subheading{Abstract} We extend Pisier's abstract version of
Grothendieck's theorem to general non-locally convex
quasi-Banach spaces.  We also prove a related result on
factoring operators through a Banach space and apply our
results to the study of vector-valued inequalities for Sidon
sets.  We also develop the local theory of (non-locally
convex) spaces with duals of weak cotype 2.

\vskip2truecm

\subheading{1.  Introduction}

In [16] (see also [18]) Pisier showed that if $X$ and $Y$
are Banach spaces so that $X^*$ and $Y$ have cotype 2 then
any approximable operator $u:X\to Y$ factors through a
Hilbert space.  This result (referred to as the abstract
version of Grothendieck's theorem in [18]) implies the usual
Grothendieck theorem by taking the special case
$X=C(\Omega)$ and $Y=L_1$ as explained in [18].

Our main result is that the abstract form of Grothendieck's
theorem is valid for quasi-Banach spaces.  To make this
precise let us say that an operator $u:X\to Y$ between two
quasi-Banach spaces is {\it strongly approximable} if it is
in the smallest subspace $\Cal A(X,Y)$ of the space $\Cal
L(X,Y)$ of all bounded operators which contains the
finite-rank operators and is closed under pointwise
convergence of bounded nets.  We define the dual $X^*$ of a
quasi-Banach space as the space under all bounded linear
functionals; this is always a Banach space.  Then suppose
$X,Y$ are quasi-Banach spaces so that $X^*$ and $Y$ have
cotype 2; we prove that if $u:X\to Y$ is strongly
approximable then $u$ factors through a Hilbert space.

Some approximability assumption is necessary even for Banach
spaces (cf.  [17]).  However, in our situation, such an
assumption is transparently required because $X^*$ could
have cotype 2 for the trivial reason that $X^*=\{0\}$; then
the only strongly approximable operator on $X$ is
identically zero.  We remark that there are many known
examples of nonlocally convex spaces $X$ with cotype 2 (e.g.
$L_p,L_p/H_p$ and the Schatten ideals $S_p$ when $p<1$
([20],[23]).  Examples of nonlocally convex spaces whose
dual have cotype 2 are less visible in nature, but in [7]
there is an example of such a space $X$ with an
unconditional basis so that $X^*\sim\ell_1.$

We also give a similar result for factorization through a
Banach space; in this case we require that $X^*$ embeds into
an $L_1-$space and that $Y$ has nontrivial cotype.  These
results are then applied to the study of Sidon sets.  We say
a quasi-Banach space $X$ is Sidon-regular if for every
compact abelian group $G$ and every Sidon subset
$E=\{\gamma_n\}_{n=1}^{\infty}$ of the dual group $\Gamma$
and for every $0<p\le \infty$ we have
$\|\sum_{k=1}^nx_k\gamma_k\|_{L_p(G,X)}\sim
\|\sum_{k=1}^n\epsilon_kx_k\|_{L_p(X)}$ where $(\epsilon_k)$
are the Rademacher functions on $[0,1]$.  It is a well-known
result of Pisier [15] that Banach spaces are Sidon-regular
but in [9] it is shown that not every quasi-Banach space is
Sidon-regular.  We show as a consequence of the above
factorization theorems that any space with nontrivial cotype
is Sidon-regular; this includes such spaces as the Schatten
ideals $S_p$ and the quotient spaces $L_p/H_p$ when $0<p<1.$

Our final section is motivated by the fact that the main
factorization theorem suggests that quasi-Banach spaces
whose duals have cotype 2 have special properties.  On an
intuitive level there is no reason to suspect that
properties of the dual space will influence the original
space very strongly in the absence of local convexity.
However we show that quasi-Banach spaces whose duals have
weak cotype 2 can be characterized internally by conditions
dual to the standard characterizations of weak cotype 2
spaces.  We show that, for example, that $X$ has a dual of
weak cotype 2 if and only if its finite-dimensional
quotients have uniformly bounded outer-volume ratios.

We refer to [11] for the essential background on
quasi-Banach spaces.  We will need the fact that every
quasi-normed space can be equivalently normed with an
$r$-norm where $0<r\le 1$ (the Aoki-Rolewicz theorem) i.e. a
quasi-norm satisfying $\|x+y\|^r\le \|x\|^r +\|y\|^r.$

Let us recall that the Banach envelope $\hat X$ of a
quasi-Banach space $X$ is defined to be the closure of
$j(X)$ where $j:X\to X^{**}$ is the canonical map (which is
not necessarily injective).  If $j$ is injective (i.e.  $X$
has a separating dual) we regard $\hat X$ as the completion
of $X$ with respect to the norm induced by the convex hull
of the closed unit ball $B_X.$ If $X$ is locally convex then
$d(X,\hat X)=\|I_X\|_{\hat X\to X}$ is equal to the minimal
Banach-Mazur distance between $X$ and a Banach space.

\newpage

\subheading{2.  The main factorization theorems}

Let us suppose that $X$ and $Y$ are $r$-Banach spaces where
$0<r\le 1.$ Suppose $u:X\to Y$ is a bounded linear operator.
We will define $\gamma_2(u)$ to be the infimum of
$\|v\|\|w\|$ over all factorizations $u=vw$ where $w:X\to H$
and $v:H\to Y$ for some Hilbert space $H.$ We define
$\delta(u)$ to be the infimum of $\|v\|\|w\|$ where $u=vw$
and $w:X\to B$ and $v:B\to Y$ for some Banach space $B.$ In
the special case when $u=I_X$ is the identity operator then
$\gamma_2(I_X)=d_X$ is the Euclidean distance of $X$ and
$\delta(I_X)=\delta_X=d(X,\hat X)$ (cf.  [5], [14]) is the
distance of $X$ to its Banach envelope.

Let $D_N=\{-1,+1\}^N$ be equipped with normalized counting
measure $\lambda$ and define the Rademacher functions
$\epsilon_i(t)=t_i$ on $D_N$ for $1\le i\le N.$ We define
$T_2^{(N)}(u)$ to be the least constant such that $$
\left(\int \|\sum_{i=1}^N
\epsilon_iu(x_i)\|^2d\lambda\right)^{1/2} \le T_2^{(N)}(u)
\left(\sum_{i=1}^N\|x_i\|^2\right)^{1/2}$$ for
$x_1,\ldots,x_N\in X.$ We define $C_2^{(N)}(u)$ to be the
least constant so that $$ \left( \sum_{i=1}^N
\|u(x_i)\|^2\right)^{1/2} \le C_2^{(N)}(u)\left(\int
\|\sum_{i=1}^N \epsilon_ix_i\|^2d\lambda\right)^{1/2} $$ for
$x_1,\ldots,x_N\in X.$ We let $T_2(u)=\sup_N T_2^{(N)}(u)$
be the type 2 constant of $u$ and $C_2(u)=\sup_N
C_2^{(N)}(u)$ be the cotype two constant of $u.$ In the case
when $u=I_X$ we let $T_2(I_X)=T_2(X)$ and $C_2(I_X)=C_2(X)$
the type two and cotype two constants of $X.$

Finally we let $K^{(N)}(u)$ be the least constant so that if
$f\in L_2(D_N,X)$ then $$ \|\sum_{i=1}^N (\int
\epsilon_iu\circ f d\lambda)\epsilon_i\|_{L_2(D_N,Y)} \le
K^{(N)}(u) \|f\|_{L_2(D_N,X)}.$$ We then let $K(u)=\sup_N
K^{(N)}(u)$.  If $u=I_X$ then $K(I_X)=K(X)$ is the
K-convexity constant of $X.$

We will need the following estimate.

\proclaim{Lemma 1} If $0<r<1$ then there is a constant
$C=C(r)$ so that for any $r$-normed space $X$, we have $
K(X)\le Cd_X^{\phi}(1+\log d_X),$ where
$\phi=(1/r-1)/(1/r-2).$ \endproclaim

\demo{Proof}It is clear that $K(X)\le \delta_XK(\hat X).$
Now we have $\delta_X\le Cd_X^{\phi}$ where $C=C(r)$.  Lemma
3 of [5].  We also have, by a result of Pisier ([16], [18])
that $K(\hat X) \le C(1+\log d_{\hat X}).$ It remains to
observe that $d_{\hat X}\le d_X$ since any operator $u:X\to
H$ where $H$ is a Hilbert space factorizes through the
Banach envelope of $X$ with preservation of
norm.\bull\enddemo

We next discuss some aspects of Lions-Peetre interpolation
(see [1] or [2] ). We will only need to interpolate between
pairs of equivalent quasi-norms on a fixed quasi-Banach
space.  Let us suppose that $X$ is an $r$-Banach space for
quasi-norm $\|\,\|_0$ and that $\|\,\|_1$ is an equivalent
$r$-norm on $X$; we write $X_j=(X,\|\,\|_j)$ for $j=0,1.$
Let $$K_s(t,x)=\inf\{(\|x_0\|_0^s+t^s\|x_1\|_1^s)^{1/s}:\
x=x_0+x_1\}$$ where $r\le s<\infty.$ Then $K_s$ is an
$r$-norm on $X.$ We define $$\|x\|_{\theta,2}=
(\theta(1-\theta))^{1/2}\left(\int_0^{\infty}\frac
{K_2(t,x)^2}{t^{1+2\theta}}dt\right)^{1/2}$$ for
$0<\theta<1.$ We write $X_{\theta,2}=(X,\|\
\|_{\theta,2})=(X_0,X_1)_{\theta,2}.$

We will need some well-known observations.

\proclaim{Lemma 2}There exists a $C=C(r)$ so that if $\|\
\|_0=\|\ \|_1$ then for any $x\in X$ we have $C^{-1}\|x\|_0
\le \|x\|_{\theta,2} \le C\|x\|_0.$\endproclaim

\demo{Proof}This follows from the simple observation that
$K_r(t,x)=\min(1,t)\|x\|_0$ and that $2^{1/2-1/r}K_r\le
K_2\le K_r.$\bull\enddemo

\proclaim{Lemma 3}There exists $C=C(r)$ so that if
$N\in\bold N$ then
$(\ell_2^N(X_0),\ell_2^N(X_1))_{\theta,2}$ is isometrically
isomorphic to $\ell_2^N(X_{\theta,2}).$\endproclaim

\demo{Proof}This follows from the routine estimate $$
K_2(t,(x_1,\ldots,x_N)) =
(\sum_{i=1}^NK_2(t,x_i)^2)^{1/2}\bull$$\enddemo

The following lemma is standard.

\proclaim{Lemma 4}Suppose $(X_0,X_1)$ are as above and that
$(Y_0,Y_1)$ is a similar pair of $r$-normings of a
quasi-Banach space $Y$.  Let $u:X\to Y$ be a bounded linear
operator.  Then $$ \|u\|_{X_{\theta,2}\to Y_{\theta,2}}\le
\|u\|_{X_0\to Y_0}^{1-\theta}\|u\|_{X_1\to Y_1}^{\theta}.$$
\endproclaim

We now combine these results to give a criterion for
convexity of the interpolated space.  For convenience we
will drop the subscript 2 and write $X_{\theta}.$ We also
recall the definition of equal norms type $p$ for $p\le 2.$
We let $\hat T_p(X)$ be the least constant so that for any
$N$ and any $x_1,\ldots,x_N\in X$ we have $$
\|\sum_{i=1}^N\epsilon_ix_i\|_{L_p(D_N,X)} \le \hat
T_p(X)N^{1/p} \max_{1\le i\le N}\|x_i\|.$$

\proclaim{Lemma 5}Suppose $(X_0,X_1)$ are as above.  Suppose
$1\le a<\infty$ and $0<\theta<r/(2-r).$ There is a constant
$C=C(a,\theta,r)$ so that if $T_2(X_0)\le a,$ then
$\delta_{X_{\theta}}\le C.$\endproclaim

\demo{Proof}Consider the map $u:\ell_2^N(X)\to
L_2(\Omega_N,X)$ defined by $u((x_i)_{i=1}^N)=
\sum_{i=1}^N\epsilon_ix_i.$ Then $\|u\|_{X_0}\le a$.  Since
$\|\ \|_1$ is an $r$-norm it follows from Holder's
inequality that $\|u\|_{X_1} \le N^{1/r-1/2}.$ Hence
$\|u\|_{X_{\theta}} \le a^{1-\theta}N^{\theta(1/r-1/2)}.$
Now $\theta(1/r-1/2)=1/2-\phi$ where $\phi>0.$ Assume
$x_1,\ldots,x_N\in X$.  Then $$ \left(\int
\|\sum_{i=1}^N\epsilon_ix_i\|_{\theta}^2d\lambda\right)^{1/2}\le
a^{1-\theta}N^{1-\phi}\max_{1\le i\le N}\|x_i\|_{\theta}.$$
This means that $\hat T_p(X_{\theta})\le a^{1-\theta}$ where
$p=(1-\phi)^{-1}>1.$ Applying Lemma 2 of [5] we get the
lemma.\bull\enddemo

We are now in position to prove the generalization of
Pisier's abstract Grothendieck theorem.

\proclaim{Theorem 6}Let $X,Y$ be quasi-Banach spaces so that
$X^*$ and $Y$ have cotype 2. Then there is a constant $C$ so
that if $u:X\to Y$ is a strongly approximable operator, then
$\gamma_2(u)\le C\|u\|.$ \endproclaim

\demo{Proof}We may suppose that both $X$ and $Y$ are
$r$-normed.  Consider first an operator $u:X\to Y$ such that
$\|u\|=1$ and $\gamma_2(u)<\infty.$ Then there is a Hilbert
space $Z$ and a factorization $u=vw$ where $w:X\to Z$ and
$v:Z\to Y$ satisfy $\|v\|\le 2\gamma_2(u)$ and $\|w\| \le
1.$

We will let $Z=Z_0$ and define $Z_1$ by the quasi-norm
$\|z\|_1=\max(\|z\|_0,\|v(z)\|_Y).$ Then $T_2(Z_0)=1$; so we
pick $0<\theta<r/(2-r)$ depending on $r$ and deduce an
estimate $\delta_{Z_{\theta}} \le C=C(r).$

Now consider the map $\tilde w_N:  L_2(D_N,X) \to
\ell_2^N(Z)$ defined by $$ \tilde w_N(f) = \left(\int w\circ
f\epsilon_i d\lambda\right)_{i=1}^N.$$ Clearly for the
Euclidean norm $\|\, \|_0$ we have $\|\tilde w_N\|_0 \le 1.$

We now consider $\|\,\|_1.$ It is routine to see that
$C_2(Z_1)\le 1+C_2(Y)\le 2C_2(Y).$ We also clearly have the
estimate $\|z\|_0\le \|z\|_1\le 2\gamma_2(u)\|z\|_0$ so that
$d_{Z_1}\le 2\gamma_2.$ From this and Lemma 1 we can obtain
an estimate $K(Z_1)\le C(\gamma_2(u))^{\phi}(1+\log
\gamma_2(u))$ where $C=C(r)$ and $\phi=(1-r)/(2-r).$ To
simplify our estimate we replace this by $K(Z_1)\le
C(\gamma_2(u))^{\tau}$ where $\tau$ depends only on $r$ and
$\phi<\tau<1.$ These estimates combine to give $$\|\tilde
w_N\|_1 \le C(\gamma_2(u))^{\tau}C_2(Y).$$

Interpolation now yields $$\|\tilde w_N\|_{\theta} \le
C(\gamma_2(u))^{\tau\theta} C_2(Y)^{\theta}.$$

Now consider $w^*:Z_{\theta}^*\to X^*.$ By taking adjoints
of $\tilde w_N$ and observing that $L_2(D_N,X)^*$ can be
identified with $L_2(D_N,X^*)$ in the standard way we see
that we have an estimate $$T_2(w^*:Z_{\theta}^*\to X^*) \le
C(\gamma_2(u))^{\tau\theta} C_2(Y)^{\theta}.$$ It follows
immediately from Maurey's extension of Kwapien's theorem
[12], [18] Theorem 3.4 that $$\gamma_2(w^*:Z_{\theta}^*\to
X^*) \le C(\gamma_2(u))^{\tau\theta}
C_2(Y)^{\theta}C_2(X^*).$$ By duality this gives the same
estimate for $\gamma_2(w:\hat X\to \hat Z_{\theta}).$ Our
previous estimate on $\delta_{Z_{\theta}}$ gives that the
norm of the identity map $I:\hat Z_{\theta} \to Z_{\theta}$
is bounded by some $C=C(r).$ Since $I_X:X\to \hat X$ has
norm one, we have:  $$\gamma_2(w:X\to Z_{\theta}) \le
C(\gamma_2(u))^{\tau\theta} C_2(Y)^{\theta}C_2(X^*).$$

Now $\|v\|_0 \le 2\gamma_2(u)$ and $\|v\|_1\le 1$ by
construction.  By interpolation we have $\|v\|_{\theta}\le
C\gamma_2(u)^{1-\theta}.$ Now by factoring through
$Z_{\theta}$ we obtain $$\gamma_2(u) \le
C(\gamma_2(u))^{1-\theta+\tau\theta}C_2(Y)^{\theta}C_2(X^*)$$
and so $$ \gamma_2(u) \le
C(C_2(Y))^{1/(1-\tau)}(C_2(X^*))^{1/(1-\tau)\theta}.$$

Thus we conclude that if $\gamma_2(u)<\infty$ then
$\gamma_2(u)\le C\|u\|$ where $C$ is a constant depending
only on $X,Y.$ The remainder of the argument is standard.
Let $\Cal J$ be the subspace of $\Cal L(X,Y)$ of all
operators for which $\gamma_2(u)<\infty.$ Then $\Cal J$
contains all finite-rank operators.  We show it is closed
under pointwise convergence of bounded nets.  Let
$(u_{\alpha})$ be a bounded net in $\Cal J$ converging
pointwise to $u.$ Then $\sup \gamma_2(u_{\alpha})=B<\infty.$
For each $\alpha$ there is a Eulcidean seminorm (i.e. a
seminorm obeying the parallelogram law) $\|\ \|_{\alpha},$
on $X$ satisfying $\|u_{\alpha}(x)\| \le \|x\|_{\alpha} \le
B\|x\|$ for $x\in X.$ By a straightforward compactness
argument there is a Euclidean seminorm $\|\ \|_E$ on $X$
satisfying $\|u(x)\|\le \|x\|_E \le B\|x\|$ for $x\in X$,
i.e.  $u\in \Cal J.$\bull\enddemo

We will next prove a similar result for factorization
through a Banach space.  We recall that a quasi-Banach space
has cotype $q$ where $q\ge 2$ if there is a constant $C$ so
that for every $N$ and all $x_1,\ldots,x_N\in X$ we have:
$$ (\sum_{i=1}^N\|x_i\|^q)^{1/q}\le C
\|\sum_{i=1}^N\epsilon_ix_i\|_{L_q(D_N,X)}.$$

\proclaim{Lemma 7}Let $X$ be a Banach space so that $X^*$ is
isomorphic to a subspace of an $L_1-$space, and suppose $Y$
is a quasi-Banach space of cotype $q<\infty.$ Then there is
a constant $C=C(X,Y)$ so that if $u:X\to Y$ is a bounded
operator then there is a Banach space $Z$ with $T_2(Z)\le C$
and a factorization $u=vw$ where $w:X\to Z$ and $v:Z\to Y$
satisfy $\|v\|\|w\| \le C\|u\|.$\endproclaim

\demo{Proof}By assumption, there is a compact Hausdorff
space $\Omega_0$ and an open mapping $q:C(\Omega)\to
X^{**}.$ It follows that there is a constant $C_0$ so that
if $E$ is a finite-dimensional subspace of $X$ there is a
finite rank operator $t_E:C(\Omega) \to X^{**}$ with
$\|t_E\|< C_0$ and if $x\in E$ with $\|x\|=1$ there exists
$f\in C(\Omega)$ with $\|f\|< 2$ and $t_Ef=x.$ It follows
from the Principle of Local Reflexivity (cf.  [22] p.76)
that we can suppose that $t_E$ has range in $X$.

We now form an ultraproduct of $X$ and $Y$.  Let $I$ be the
the collection of all finite-dimensional subspaces $E$ of
$X$ and let $\Cal U$ be an ultrafilter on $I$ containing all
sets of the form $\{E:E\supset F\}$ where $F$ is a fixed
finite-dimensional subspace.  Consider the space $X_{\Cal
U}$ defined to be the quotient of $\ell_{\infty}(I;X)$ by
the subspace $c_{\Cal U,0}(I;X)$ of all $(x_E)$ such that
$\lim_{\Cal U}x_E=0;$ $X_{\Cal U}$ is thus the space of
(equivalence classes of) $(x_E)$ normed by $\lim_{\Cal
U}\|x_E\|.$ We regard $X$ as a subspace of $X_{\Cal U}$ by
identifying $x$ with the constant function $x_E=x$ for all
$E.$ We similarly introduce $Y_{\Cal U}$ and note that
$Y_{\Cal U}$ has cotype $q$ with the same constant as $Y$.
We extend $u:X\to Y$ to $u_1:X_{\Cal U}\to Y_{\Cal U}$ by
setting $u_1 ((x_E)_{E\in I})=(u(x_E)_{E\in I}).$ Let us
also introduce the operator $t:C(\Omega)\to X_{\Cal U}$ by
putting $t(f)=(t_E(f))_{E\in I}.$ Consider
$u_1t:C(\Omega)\to Y_{\Cal U}.$ By Theorem 4.1 of [10] there
is a regular probabilty measure $\mu$ on $\Omega$ so that if
$p=q+1$ then $\|u_1t(f)\|\le C\|u\|(\int |f|^pd\mu)^{1/p}.$
Here $C$ depends on $X,Y$ but not on $u.$ This implies that
$u_1t=v_1j$ where $j:C(\Omega)\to L_p(\mu)$ is the canonical
injection and $\|v_1\|\le C\|u\|.$ Let $N=v_1^{-1}(0)$ and
form the quotient $Z_1=L_p/N;$ let $\pi$ be the quotient
map.  For each $x\in X$ there exists $f\in C(\Omega)$ and
$t(f)=x;$ this follows from the choice of ultrafilter.  Then
$w(x)=\pi j(f)$ is uniquely determined independent of $f.$
Furthermore $f$ can be chosen so that $\|f\|\le 2\|x\|$ so
that $\|w\|\le 2.$ Let $Z$ be the closure of the range of
$w$.  Then clearly $v_1(\pi^{-1}(Z))\subset Y$ so that we
can define $v:Z\to Y$ with $\|v\|\le C$ and $vw=u.$ Finally
$T_2(Z)\le T_2(Z_1)\le T_2(L_p)$ is bounded by a constant
depending only on $q.$\bull\enddemo

\proclaim{Theorem 8} Suppose $X$ is a quasi-Banach space
such that $X^*$ is isomorphic to a subspace of an
$L_1-$space, and $Y$ is a quasi-Banach space of cotype
$q<\infty.$ There is a constant $C$ so that if $u:X\to Y$ is
a strongly approximable operator then $\delta(u)\le C\|u\|.$
\endproclaim

\demo{Proof}Let us assume that $X,Y$ are both $r$-normed.
Suppose first that $u:X\to Y$ satisfies $\|u\|=1$ and
$\delta(u)<\infty.$ We show that $\delta(u)\le C$ where $C$
depends only on $X,Y$.  By Lemma 7, $u$ can be factored
through a Banach space $Z$ satisfying $T_2(Z)\le C$ where
$C=C(X,Y)$ so that $u=vw$ where $w:X\to Z$ with $\|w\|=1$
and $v:Z\to Y$ with $\|v\|\le C\delta(u).$ Let $Z=Z_0$ and
introduce $Z_1$ by setting $\|z\|_1=\max(\|z\|,\|v(z)\|).$
If we pick $\theta<r/(2-r)$ then $\delta_{Z_{\theta}}\le C$
where again $C$ depends only on $X,Y.$ Hence $$\delta(u) \le
C\|v\|_{Z_{\theta}\to Y}\|w\|_{X\to Z_{\theta}}.$$ By
interpolation this yields $$\delta(u) \le
C(\delta(u))^{1-\theta},$$ and hence $\delta(u)\le C.$ The
remainder of the proof is similar to that of Theorem
6.\bull\enddemo

\vskip2truecm

\subheading{3.  Applications to Banach envelopes and Sidon
sets}

It is proved in [7] that if $X$ is a natural quasi-Banach
space (i.e. a space isomorphic to a subspace of a space
$\ell_{\infty}(I;L_p(\mu_i))$ with the strong approximation
property and if $Y$ is any subspace of $X$ such that $Y^*$
has cotype $q<\infty$ then $Y$ is locally convex.  We
present now two variations on this theme.

Let us say that a quasi-Banach space $X$ is (isometrically)
subordinate to a quasi-Banach space $Y$ if $X$ is
(isometrically) isomorphic to a closed subspace of a space
$\ell_{\infty}(I;Y)$ for some index set $I.$ Thus a
separable space $X$ is natural if it is subordinate to
$L_p[0,1]$ for some $0<p<1.$

\proclaim{Theorem 9}Let $Z$ be a quasi-Banach space and let
$X$ be subordinate to $Z.$ Assume that either $X$ or $Z$ has
the strong approximation property.  Let $Y$ be any subspace
of $X.$ Then\newline (1) If $Z$ has cotype 2 and if $Y^*$
has cotype 2 then $Y$ is locally convex.\newline (2) If $Z$
has cotype $q<\infty$ and $Y^*$ is isomorphic to a subspace
of an $L_1-$space, then $Y$ is locally convex.  \endproclaim

\demo{Proof}The proofs are essentially identical.  We
therefore prove only (2).  Let $j:Y\to \ell_{\infty}(Z)$ be
the inclusion map.  Then since $j$ factors through $X$ it is
strongly approximable, under either hypothesis.  Let
$\pi_i:\ell_{\infty}(Z)\to Z$ be the co-ordinate map.  Then
by Theorem 8, we have $\delta(\pi_ij)\le C$ for some
constant $C$ depending only on $Y.$ Thus for $y\in Y,$
$\|y\|=\sup_i\|\pi_ijy\| \le C\|y\|_{\hat Y}.\bull$\enddemo

Let us now give an application.  Suppose $G$ is a compact
abelian group with normalized Haar measure $\mu_G,$ and
suppose $\Gamma$ is the dual group.  We recall that a subset
$E\subset\Gamma$ is a Sidon set if for every
$\epsilon_{\gamma}=\pm 1$ there exists $\nu\in \Cal M(G)$
whose Fourier transform satisfies $\hat
\nu(\gamma)=\epsilon_{\gamma}.$

In [9] the first author introduced the property $\Cal
C_p(X)$ for a subset $E$ of $\Gamma$ where $0<p\le \infty.$
We say that $E$ has $\Cal C_p(X)$ if there is a constant $M$
so that for any $\gamma_1,\ldots,\gamma_n\in E$ and any
$x_1,\ldots,x_n\in X$ we have $$
M^{-1}\|\sum_{k=1}^nx_k\epsilon_k\|_{L_p(D_n,X)} \le
\|\sum_{k=1}^nx_k\gamma_k\|_{L_p(G,X)} \le
M\|\sum_{k=1}^nx_k\epsilon_k\|_{L_p(D_n,X)}$$ where
$\epsilon_1,\ldots,\epsilon_n$ are the Rademacher functions
on $D_n.$ Let us say that a quasi-Banach space $X$ is {\it
Sidon-regular} if every Sidon set $E$ has property $\Cal
C_p(X)$ for every $0<p\le \infty.$ It is a well-known result
of Pisier [15] that every Banach space is Sidon-regular.  By
way of contrast, in [9] an example of a quasi-Banach space
which is not Sidon-regular is constructed.  However, every
natural space is Sidon-regular.  The above results enable us
to extend this to a wider class of spaces.

\proclaim{Theorem 10}Let $X$ be a quasi-Banach space of
cotype $q<\infty.$ Then every quasi-Banach space which is
subordinate to $X,$ (and, in particular, $X$ itself) is
Sidon-regular.  \endproclaim

\demo{Proof}This is very similar to the proof of Theorem 4
in [9].  Suppose $G$ is a compact Abelian group and $E$ is a
Sidon subset of $\Gamma.$ Let $s=\min(p,2).$ Let
$Z=L_s(G,X).$ Then $Z$ also has cotype $q$.  To see this we
need first to observe that the Kahane-Khintchine inequality
holds in an arbitrary quasi-Banach space (Theorem 2.1 of
[6]) so that there is a constant $C$ depending only on $X$
so that if $x_1,\ldots,x_n\in X$ then
$$\left(\int\|\sum_{k=1}^n\epsilon_kx_k\|^qd\lambda\right)^{1/q}
\le
C\left(\int\|\sum_{k=1}^n\epsilon_kx_k\|^sd\lambda\right)^{1/s}.$$
Now if $f_1,\ldots,f_n\in Z$ then, for constants
$C_1,C_2,C_3$ depending only on $X,$ $$ \align
\left(\sum_{k=1}^n\|f_k\|^q\right)^{1/q} &\le
C_1\left(\int_G
(\sum_{k=1}^n\|f_k(t)\|_X^q)^{s/q}d\mu_G(t)\right)^{1/s}\\
&\le
C_2\left(\int_G\left(\int_{D_n}\|\sum_{i=1}^n\epsilon_kf_k(t)\|^q_X
d\lambda\right)^{s/q} d\mu_G(t)\right)^{1/s}\\ &\le
C_3\left(\int_{D_n}\int_G
\|\sum_{k=1}^n\epsilon_kf_k(t)\|^s_X
d\mu_G(t)d\lambda\right)^{1/s}\\ &\le C_3
\|\sum_{k=1}^n\epsilon_kf_k\|_{L_q(D_n,Z)}.  \endalign $$

Now let $E_n$ be any sequence of finite subsets of $E.$ Let
$\Cal P_{E_n}(Y)$ be the space of $Y$-valued polynomials
$\sum_{\gamma\in E_n}y_{\gamma}\gamma$ equipped with the
$L_p(G,Y)$ quasi-norm.  Then $\Cal P_{E_n}(Y)$ is
isometrically subordinate to $Z.$ Next equip the
finite-dimensional space $\ell_{\infty}(E_n)$ of all bounded
functions $h:E_n\to\bold C$ with the quasi-norm of the
operator $T_h:\Cal P_{E_n}\to\Cal P_{E_n}$ given by
$T_h(\sum y_{\gamma}\gamma)=\sum h(\gamma)y_{\gamma}\gamma.$
Let us denote this space $\Cal M_n.$ Then $\Cal M_n$ is
isometrically subordinate to $Z.$ Thus the product $c_0(\Cal
M_n)$ is isometrically subordinate to $Z$ and has the strong
approximation property.  However, as in Theorem 4 of [9] the
assumption that $E$ is a Sidon set shows that we have a
constant $C$ depending only the Sidon constant of $E$ so
that the envelope norm on $\Cal M_n$ satisfies
$\|h\|_{\infty}\le \|h\|_{\hat{\Cal M}_n} \le
C\|h\|_{\infty}.$ Thus the envelope of $c_0(\Cal M_n)$ is
isomorphic to $c_0$ and Theorem 9(ii) applies to give that
this space is locally convex so that for some uniform
constant $C'$ we have for every $n,$ $\|h\|_{\Cal M_n} \le
C'\|h\|_{\infty}.$ As in [9] Theorem 4 this implies that $E$
has property $\Cal C_p(Y).$\bull\enddemo

\demo{Remarks}The above theorem applies to $L_p/H_p$ when
$p<1$ and to the Schatten ideals $S_p$ when $p<1,$ since
these spaces have cotype 2 by recent results of Pisier [20]
and Xu [23].  These spaces are known not to be natural;
$S_p$ is A-convex (i.e. has an equivalent plurisubharmonic
quasi-norm) while $L_p/H_p$ is not A-convex (see [8]).

Let us also remark that if $0<p<1$ and $E$ is a symmetric
$p$-convex sequence space with the Fatou property then we
can define an associated Schatten class $S_E$ (see
Gohberg-Krein [4] for the Banach space versions).  Precisely
if $H$ is a separable Hilbert space and $A$ is a compact
operator with singular values $(s_n(A))$ we say $A\in S_E$
if $(s_n(A))\in E$ and we set $\|A\|_E=\|(s_n(A))\|_E.$ It
can then be shown that $S_E$ is subordinate to $S_p$.  In
fact we define a sequence space $F$ by $\|(t_n)\|_F = \sup
\{\|(s_nt_n)\|_p:\|(s_n)\|_E\le 1\}$ and it can then be
shown that $\|A\|_E =\sup\{ \|AB\|_p:\|B\|_F\le 1\}.$ This
result follows quickly from an inequality of Horn (cf.  [4]
pp. 48-9) that $$ \sum_{j=1}^ks_j(AB)^p \le
\sum_{j=1}^ks_j(A)^ps_j(B)^p$$ for every $k.$

\vskip2truecm

\subheading{4.  Quasi-Banach spaces with duals of weak
cotype 2}

Let $X$ be a finite-dimensional continuously quasi-normed
quasi-Banach space with unit ball $B_X$.  We recall that the
{\it volume-ratio} of $X$ is defined by $\vr(X)=(\vol
B_X/\vol\Cal E)^{1/n}$ where $\Cal E$ is an ellipsoid of
maximal volume contained in $B_X$ and $n=\dim X.$ We define
the {\it outer volume-ratio} of $X$ by $\ovr(X)=(\vol\Cal
F/\vol B_X)^{1/n}$ where $\Cal F$ is an ellipsoid of minimal
volume contaning $B_X.$ The Santalo inequality ([19], [21])
shows that $\ovr(X)\ge (\vol\Cal B_{X^*}/\vol
F^0)^{1/n}=\vr(X^*).$ The reverse Santalo inequality of
Bourgain and Milman ([3],[19]) shows that, if $X$ is normed,
$\ovr(X)\le C\vr(X^*)$ so that $\ovr(X)$ is then equivalent
to $\vr(X^*).$ For general quasi-normed spaces the reverse
Santalo inequality is not available.

We recall that a Banach space $X$ is of weak cotype 2 if
there exists $C$ so that whenever $H$ is a
finite-dimensional Hilbert space with orthonormal basis
$(e_1,\ldots,e_n)$ and $u:H\to X$ is a linear operator then
$a_k(u)\le Ck^{-1/2}\ell(u)$ for $1\le k\le n.$ Here $$
\ell(u)=(\bold E(\|\sum_{k=1}^ng_ku(e_k)\|^2))^{1/2}$$ (for
$g_1,\ldots,g_n$ a sequence of independent normalized
Gaussian random variables) and $a_k(u)=\inf\{\|u-v\|:v:H\to
X,\ \text{rank }v<k\}.$ The least such constant $C$ is
denoted by $wC_2(X).$ It is known that $X$ is of weak cotype
2 if and only if there exists $C$ so that $\vr(E)\le C$ for
every finite-dimensional subspace of $X.$ See Pisier [19]
for details.  It follows quickly that $X^*$ is of weak
cotype 2 if and only if $\ovr(E)$ is bounded for all
finite-dimensional quotients of $X.$ We prove in this
section that the same characterization extends to
quasi-Banach spaces.

We will require a preparatory lemma:

\proclaim{Lemma 11}Let $E$ be an $N$-dimensional Euclidean
space and suppose $B$ is the unit ball of an $r$-norm on
$E$.  Let $S$ be a subspace of $E$ of dimension $k.$ Suppose
$1/r=\beta\in\bold N.$ Then $$ \frac{\vol (B\cap S)\vol
P_{S^{\perp}}(B)}{\vol B}\le \binom{N\beta}{k\beta}$$ where
$P_{S^{\perp}}$ is the orthogonal projection of $E$ onto
$S^{\perp}.$\endproclaim

\demo{Proof}We duplicate the argument of Lemma 8.8 of [19]
(p. 132).  One finds that $\vol B\ge a\vol (B\cap S)\vol
P_{S^{\perp}}B$ where $$ \align a&= (N-k)\int_0^1
(1-t^r)^{k/r}t^{N-k-1}dt\\ &= \frac{N-k}r \int_0^1
(1-s)^{k/r}s^{N/r-k/r-1}ds\\ &=\binom{N\beta}{k\beta}.\bull
\endalign $$ \enddemo

\proclaim{Lemma 12}There is a constant $C$ depending only
$r$ so that if $E$ is a finite-dimensional $r$-normed space
then $d_E \le Cd_{\hat E}^ {2/r-1},$ and $\delta_E \le
Cd_{\hat E}^{2/r-2}.$ \endproclaim \demo{Proof}By Lemma 3 of
[5] we have $d_E \le \delta_Ed_{\hat E}\le
Cd_E^{\phi}d_{\hat E}$ where $\phi=(1/r-1)/(1/r-1/2).$ This
proves the first part and the second part follows on
reapplying Lemma 3 of [5].  \bull\enddemo

\proclaim{Theorem 13}Let $X$ be a quasi-Banach space and
suppose $0<\alpha<1.$ Then $X^*$ has weak cotype 2 if and
only if there is a constant $C$ so that whenever $F$ is a
finite-dimensional quotient of $X$, there exists a quotient
$E$ of $F$ with $\dim E\ge \alpha\dim F$ and $d_E\le
C.$\endproclaim

\demo{Proof}Suppose $X^*$ has weak cotype 2. Then if $F$ is
a finite-dimensional quotient, $wC_2(F^*)\le wC_2(X^*)$ and
so has a subspace $G$ with $\dim G\ge \alpha\dim F$ and
$d_G\le C$ (where $C$ depends only on $X$ and $\alpha$).
Let $E=F/G^{\perp}.$ Then $d_{\hat E}=d_G$ and so the
preceding Lemma gives an estimate $d_E\le C'$ where
$C'=C'(\alpha,X).$

Conversely if $X$ has the given property then it is easy to
see that $\hat X$ must also have the same property and this
leads quickly to the fact that $X^*$ has weak cotype 2 by
the results of [13].\bull\enddemo

\proclaim{Proposition 14}Suppose $0<r<1$ and $a\ge 1$; then
there is a constant $C=C(a,r)$ so that if $E$ is an
$N$-dimensional $r$-normed space and $wC_2(E^*)\le a$ then
$\ovr(E)\le C.$ \endproclaim

\demo{Proof}In the argument which follows we use $C$ for a
constant which depends only $a,r$ but may vary from line to
line.  It suffices to establish the result when
$1/r=\beta\in\bold N.$ Let $\Cal E$ be the ellipsoid of
minimal volume containing $B_E$.  Using this ellipsoid to
introduce an inner-product we can define
$\|x\|_{E^*}=\sup_{b\in B_E} |(x,b)|$.  Then $\Cal E$ is the
ellipsoid of minimal volume containing $B_{\hat E}=\text{co
}B_E$ and the ellipsoid of maximal volume contained in
$B_{E^*}.$

Now, by imitating the argument of Theorem 8 of [5] we can
construct an increasing sequence of subspaces
$(W_k)_{k=1}^{\infty}$ of $E$ with $\dim W_k=N-\sigma_k\ge
(1-2^{-k})N$ and $B_{E^*}\cap W_k\subset C2^{3k}\Cal E$ for
$k\ge 1.$ We let $\tau_k=\sigma_k-\sigma_{k+1}.$

It follows from the Hahn-Banach theorem that $\Cal E\supset
P_{W_k}B_{\hat E}\supset C^{-1}2^{-3k}\Cal E\cap W_k.$ Now,
identifying $H_k=E/W_k^{\perp}$ with $W_k$ under the
quasinorm with unit ball $P_{W_k}(B_E)$ this implies that
$d_{\hat H_k}\le C2^{3k}.$ Now from Lemma 12
$\delta_{H_k}\le C2^{sk}$ for suitable $s>0$ depending on
$r.$ We conclude that $\Cal E\cap W_k \subset C2^{t
k}P_{W_k}(B_{E}),$ where $t$ depends only on $r,$ and
$C=C(a,r).$

Let $Z_k$ is the orthogonal complement of $W_k$ in
$W_{k+1}.$ Notice that $\dim Z_k=\tau_k.$ Now by Lemma 11,
if we set $A_k=P_{W_{k+1}}(B_E)\cap Z_k,$ $$ \vol
P_{W_k}(B_E)\vol A_k\le
\binom{(N-\sigma_{k+1})\beta}{\tau_k\beta} \vol
P_{W_{k+1}}(B_E).$$

Let $l$ be the first index for which $W_l=E$ and so
$\sigma_l=0.$ We first estimate $$ \align
\prod_{k=1}^{l-1}\binom{(N-\sigma_{k+1})\beta}{\tau_k\beta}&=
\frac{(N\beta)!}{((N-\sigma_1)\beta)!(\tau_1\beta)!\ldots
(\tau_{l-1}\beta)!}\\
&=\binom{N\beta}{\sigma_1\beta}\prod_{k=1}^{l-2}\binom{\sigma_k\beta}
{\tau_k\beta}\\ &\le 2^{\beta(N+\sum\sigma_k)}\\ &\le
2^{2\beta N} \endalign $$

Now $\Cal E\cap Z_k \subset C2^{t(k+1)}A_k$ and so
$\log_2\vol P_{Z_k}(B_E)\ge -Ck\tau_k+\log_2\vol \Cal E\cap
Z_k.$ (If $\tau_k=0$ we interpret the relative volume as
one).  Summing we obtain since $\tau_k\le 2^{-k}N,$ $$
\sum_{k=1}^{l-1}\log_2\vol A_k \ge
-CN+\sum_{k=1}^{l-1}\log_2\vol \Cal E\cap Z_k.$$

Thus $$ \align \log_2\vol B_E &= \log_2\vol P_{W_l}(B_E)\\
&\ge \log_2\vol P_{W_1}(B_E) - CN + \sum_{k=1}^{l-1}\vol
\Cal E\cap Z_k\\ &\ge -CN +\log_2\vol \Cal E\cap W_1
+\log_2\vol\Cal E\cap W_1^{\perp}\\ &\ge -CN +\log_2\vol\Cal
E. \endalign $$ This completes the proof of the Proposition.

\demo{Remark}This Proposition can be interpreted as follows.
Suppose $X$ is a finite-dimension\-al normed space so that
$wC_2(X^*)\le a$.  Consider the set $\partial B_X$ of
extreme points of $B_X$ and form the $r$-convex hull
$A_r=\text{co}_r\partial B_X.$ Then although $A_r$ is
smaller than $B_X$ it is not too much smaller, for $\vol
B_X/\vol A_r\le C^{\dim E}.$\enddemo

\proclaim{Theorem 15}Let $X$ be a quasi-Banach space.  Then
$X^*$ has weak cotype 2 if and only there is a constant $C$
so that $\ovr(E)\le C$ for every finite-dimensional quotient
$E$ of $X.$\endproclaim

\demo{Proof}First suppose $\ovr(E)\le C$ for every
finite-dimensional quotient $E$ of $X.$ Let $F$ be a
finite-dimensional subspace of $X^*$ and consider
$E=X/F^{\perp}.$ It is easy to see that the envelope norm on
$E$ is the quotient norm from the envelope norm on $X.$
Clearly from the definition, $\ovr(\hat E)\le\ovr(E)\le C.$
Hence by the Santalo inequality $\vr(F)\le C$.  This shows
that $X^*$ has weak cotype 2.

The converse is immediate from Proposition 14.\bull\enddemo

\vskip2truecm

\subheading{References}

\item{1.}C.  Bennett and R. Sharpley, {\it Interpolation of
operators,} Academic Press, Orlando 1988.

\item{2.}J.  Bergh and J. L\"ofstr\"om, {\it Interpolation
spaces.  An Introduction,} Springer, New York 1976.

\item{3.}J.  Bourgain and V.D.  Milman, New volume ratio
properties for convex symmetric bodies in $\bold R^n,$
Inventiones Math. 88 (1987) 319-340.

\item{4.}I.C.  Gohberg and M.G. Krein, {\it Introduction to
the theory of linear non-selfadjoint operators,} Amer.
Math.  Soc.  Trans. 19, 1969.

\item{5.}Y.  Gordon and N.J.  Kalton, Local structure theory
for quasi-normed spaces, Bull.  Sci.  Math. to appear.

\item{6.}N.J.  Kalton, Convexity, type and the three space
problem, Studia Math. 69 (1981) 247-287.

\item{7.}N.J.  Kalton, Banach envelopes of non-locally
convex spaces, Canad.  J. Math. 38 (1986) 65-86.

\item{8.}N.J.  Kalton, Plurisubharmonic functions on
quasi-Banach spaces, Studia Math. 84 (1986) 297-324.

\item{9.}N.J.  Kalton, On vector-valued inequalities for
Sidon sets and sets of interpolation, Colloq.  Math. to
appear

\item{10.}N.J.  Kalton and S.J.  Montgomery-Smith, Set
functions and factorization, to appear

\item{11.}N.J.  Kalton, N.T.  Peck, and J.W.  Roberts, {\it
An F-space sampler,} London Math.  Soc.  Lecture Notes 89,
Cambridge University Press 1985.

\item{12.}S.  Kwapien, Isomorphic characterization of inner
product spaces by orthogonal series with vector
coefficients, Studia Math. 44 (1972) 583-595.

\item{13.}V.D.  Milman and G. Pisier, Banach spaces with a
weak cotype 2 property, Israel J. Math., 54 (1986) 139-158.

\item{14.}N.T.  Peck, Banach-Mazur distances and projections
on $p$-convex spaces, Math.  Zeit. 177 (1981) 131-142.

\item{15.}G.  Pisier, Les in\'egalit\'es de Kahane-Khintchin
d'apr\`es C. Borell, S\'eminaire sur la g\'eometrie des
\'espaces de Banach, Ecole Polytechnique, Palaiseau,
Expos\'e VII, 1977-78.

\item{16.}G.  Pisier, Un th\'eor\`eme sur les op\'erateurs
entre espaces de Banach qui se factorisent par un espace de
Hilbert, Ann.  \'Ecole Norm.  Sup. 13 (1980) 23-43.

\item{17.}G.  Pisier, Counterexamples to a conjecture of
Grothendieck, Acta Math. 151 (1983) 181-208.

\item{18.}G.  Pisier, {\it Factorization of linear operators
and geometry of Banach spaces,} NSF-CBMS Regional Conference
Series, No. 60, American Mathematical Society, Providence,
1986.

\item{19.}G.  Pisier, {\it The volume of convex bodies and
geometry of Banach spaces,} Cambridge Tracts 94, Cambridge
University Press, 1989.

\item{20.}G.  Pisier, A simple proof of a theorem of Jean
Bourgain, to appear

\item{21.}L.  Santal\'o, Un invariante afin para los cuerpos
convexos del espacio de $n$ dimensiones, Portugal Math. 8
(1949) 155-161.

\item{22.}P.  Wojtaszczyk, {\it Banach spaces for analysts,}
Cambridge Studies 25, Cambridge University Press 1991.

\item{23.}Q.  Xu, Applications du th\'eor\`emes de
factorisation pur des fonctions \`a valeurs operateurs,
Studia Math. 95 (1990) 273-292.

\bye